\documentclass[12pt]{article}
\pdfoutput=1
\usepackage[ascii]{inputenc}
\usepackage[T1]{fontenc}
\usepackage[english]{babel}
\usepackage{amsmath,amssymb,amsfonts,textcomp,epsfig}
\usepackage{hyperref}
\usepackage{color}
\usepackage{calc}
 \usepackage[all]{xy}
\def\Im{ \textrm{Im}\,}
\def\R{\mathbb{R}}
\def\C{\mathbb C}

\def\P{\mathbb P}
\title{On the number of limit cycles which appear by perturbation of
Hamiltonian two-saddle cycles of planar vector fields}
 \author{Lubomir Gavrilov \\
 \normalsize \it Institut de Math\'{e}matiques de Toulouse, UMR 5219\\
 \it Universit\'{e}  de Toulouse,  31062 Toulouse,  France  }
\begin{document}
\maketitle
\newtheorem{definition}{Definition}
\newtheorem{remark}{Remark}
\newtheorem{theorem}{Theorem} 
\newtheorem{lemma}{Lemma}
\newtheorem{proposition}{Proposition}
\newtheorem{corollary}{Corollary}
\vspace{5mm} \noindent 2000 MSC scheme numbers: 34C07, 34C08, 34C05
\begin{abstract}
We find an upper bound to the maximal number of limit cycles, which bifurcate
from a hamiltonian two-saddle loop of an analytic vector field, under an
analytic deformation.
\end{abstract}
\newpage
\tableofcontents

\section{Introduction}

Consider a $N$-parameter analytic family of analytic plane vector fields
$X_\lambda$, $\lambda \in(\mathbb{R}^N,0)$,  such that $X_0$ has a $k$-saddle
cycle (a hyperbolic $k$-graphic) $\Gamma_k$, as on fig.\ref{cycles}. The
cyclicity $Cycl(\Gamma_k,X_\lambda)$ of $\Gamma_k$ is, roughly speaking, the
maximal number of limit cycles of $X_\lambda$ which tend to $\Gamma_k$ as
$\lambda \rightarrow 0$.
The first results on  the cyclicity of  one-saddle connection (also called
homoclinic connection) go back to Andronov and Leontovich in 1937 (but they
were published only in 1959 \cite{anle59}).
The cyclicity $Cycl(\Gamma_1)$ has been studied later in full generality by
Roussarie \cite{rous86,rous89}, see also \cite{rous98,ilya02} for an extensive
list of references.

 The main technical tool of the Roussarie's method is an asymptotic
expansion of the Dulac map (transport map near the saddle point) $$x\rightarrow
d_\lambda(x)$$ in terms of $x^k, x^k\omega(x,\varepsilon)$ where
$\omega(x,\varepsilon)$ is the so called Ecalle-Roussarie compensator
$$
\omega(x,\varepsilon) = \frac{x^{-\varepsilon} - 1}{\varepsilon}, \omega(x,0)=
- \ln x ,   x\in (\mathbb{R}^+,0)
$$
and $\varepsilon=\varepsilon(\lambda)$ is the trace of the vector field
$X_\lambda$ at the saddle point. Let $P_\lambda$ be the Poincar\'{e} first return
map, associated to $\Gamma_1$ and $X_\lambda$.  The most delicate case to be
studied is when $\Gamma_1$ is of infinite co-dimension ($P_0=id$). For
$\lambda\sim 0$ the
 map $P_\lambda$ is composed by a Dulac map (near the saddle point) and
an analytic map (the transport map along the homoclinic orbit). The usual
derivation-division algorithm then provides an upper bound for the cyclicity in
terms of the number of the coefficients of the asymptotic series of the
displacement map $P_\lambda-id$, which vanish as $\lambda=0$. The same method
was applied more recently to one-parameter deformations of Hamiltonian
two-saddle loops (called  also heteroclinic Hamiltonian connections), under the
non-generic assumption that one of the separatrices of $\Gamma_2$ remains
unbroken \cite{duro06,cdr07}. Recall that a $k$-saddle cycle $\Gamma_k$ is said
to be Hamiltonian, provided that there is a neighborhood of $\Gamma_k$ in which
$X_0$ allows an analytic first integral with only Morse critical points.

\emph{The purpose of the present paper is to extend these results to the case
of an arbitrary analytic perturbation of a Hamiltonian two-loop $\Gamma_2$,
having two hyperbolic equilibrium points.}

Our approach is different, as we do not use the asymptotic series of the
corresponding Dulac maps $d^1_\lambda$ and $d^2_\lambda$, shown on
fig.\ref{fig4}). Recall that in the one-parameter case $\lambda=\varepsilon\in
(\mathbb{R},0)$, the displacement function $ d^1_\varepsilon-d^2_\varepsilon$
 can be approximated  by an appropriate Abelian integral $I(.)$ (or more
generally, an iterated path integral) depending on a parameter $t$ as follows
$$
d^1_\varepsilon(t)-d^2(t)_\varepsilon = \varepsilon^d I(t) + \dots,
\varepsilon\sim 0 .
$$
Therefore to count the zeros of $ d^1_\varepsilon-d^2_\varepsilon$
(corresponding to limit cycles) it is enough to count the zeros of $I(.)$,
which can be done  by making use of the so called "Petrov trick" (based on the
argument principle), see \cite{petr90} and the references given there.

The above considerations hold true at least far from the singular points of the
vector field $X_\lambda$. As discovered in \cite{duro06,cdr07}, however, not
all limit cycles in a neighborhood of a two-saddle loop can be approximated in
such a way. The missing "alien" limit cycles are moreover non-avoidable in
generic $N$-parameter deformations $X_\lambda$ with $N\geq 4$. For this reason,
we apply the argument principle directly to the displacement map $
d^1_\varepsilon(.)-d^2(.)_\varepsilon$ in an appropriate complex domain, in
order to obtain an estimate to the number of its complex zeros (corresponding
to complex limit cycles).

The main technical result of the paper is Lemma \ref{mainlemma}, which claims
that the zero locus of the imaginary part of the Dulac map is a real analytic
curve of $\mathbb{R}^2=\mathbb{C}$ at the origin. This makes possible to
investigate the number of the zeros of the imaginary part of
$d^1_\lambda-d^2_\lambda$ along the zero locus of the imaginary part of
$d^1_\lambda$. Indeed,  the intersection numbers of two analytic curves is
easily computed. The proof of our main result, Theorem \ref{main}, is then
completed by making use of the Petrov trick.

In the course of the proof of Theorem \ref{main} we assume, for the sake of
simplicity, that our deformations depend on a single small parameter
$\varepsilon$. General multi-parameter deformations $\lambda\rightarrow
X_\lambda$ of $X_0$ are then studied along the same lines, as it follows from
the Hironaka's desingularization theorem. We explain this in Appendix \ref{b2},
see Theorem \ref{main+}.

Deformations of an arbitrary (possibly non-Hamiltonian) two-loop of infinite
co-dimension ($P_\lambda=id$) can be studied in a similar way, and will be
considered in another paper.
\section{Description of the result}
\label{description}
\begin{figure}
\begin{center}
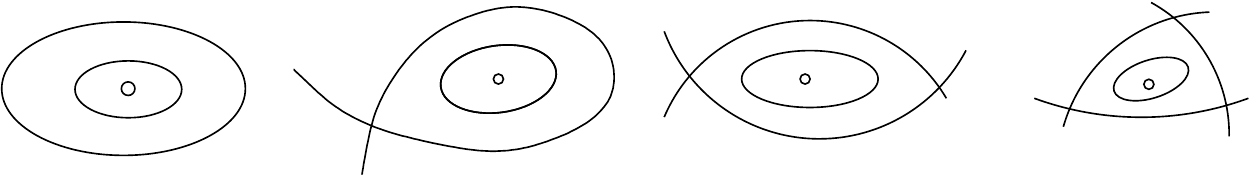
\end{center} \caption{Hamiltonian $k$-saddle cycles}\label{cycles}
\end{figure}
 Let $X_0$ be a real plane vector field. Recall that a polycycle
of $X_0$ is a topological polygon composed of separatrices and singular points.
A $k$-saddle cycle of $X_0$ (or a hyperbolic $k$-graphic) denoted $\Gamma_k$,
is a polycycle composed of $k$ distinct saddle-type singular points
$p_1,p_2,\dots,p_k$, $p_{k+1} =p_1$ and separatrices (heteroclinic  orbits)
connecting $p_i$ to $p_{i+1}$. Let $\sigma$ be a segment transversal to the
polycycle. The $k$-saddle cycle is said to be Hamiltonian, provided that $X_0$
has an analytic first integral $f$ having Morse critical points at $p_i$. It
follows that $\Gamma_k$ bounds an annulus of periodic orbits $\{(x,y):
f(x,y)=t\}_t$ of $X_0$. Thus, a Hamiltonian $0$-saddle cycle is simply a
center, a Hamiltonian $1$-saddle cycle is a homoclinic loop bounding a period
annulus, a Hamiltonian $2$-saddle cycle is a double heteroclinic loop bounding
a period annulus etc., see fig.\ref{cycles}.

One can find a "tubular neighborhood"  $U\subset \C^2\cong \R^4$ of
$$\Gamma_k\subset \{(x,y)\in\R^2 : f(x,y)=0\}$$ such that
\begin{itemize}
    \item $\bar{U}$ is compact smooth manifold with a (real three-dimensional) border.
    \item $f$ is analytic in some neighborhood of $\bar{U}$
    \item the border $\partial \bar{U}$ is transversal to the
  complex curves $\{(x,y)\in\C^2: f(x,y)=t\}$, provided that $|t|$ is
  sufficiently small.
  \item the intersection of $U$ with the singular fiber $ \{(x,y)\in\C^2 : f(x,y)=0\}$ is a union of $k$
  Riemann surfaces $D_i$, each of them homeomorphic to an open disc. $D_i$ intersects transversally $D_{i+1}$ at $p_i$,
  $i= 1\dots k$, and $D_i\cap
D_j=\emptyset$ for $|i-j|\neq 1$.
\end{itemize}
It follows that $$f: U \rightarrow \C$$ defines a locally trivial fibration
over a punctured neighborhood of the origin in $\C$, and each fiber
\begin{equation}\label{fiber}
 F_t =
U\cap\{(x,y)\in\C^2: f(x,y)=t\} ,  t\neq 0
\end{equation} is homeomorphic to a genus one
surface with $k$ punctures.

A one-parameter analytic deformation of $X_0$ is a a family $X_\varepsilon$ of
real-analytic plane vector fields, depending analytically on a real parameter
$\varepsilon\in(\mathbb{R},0)$, and defined in a suitable neighborhood of the
$k$-saddle cycle $\Gamma_k$.
 The corresponding foliation
$\mathcal{F}_\varepsilon$ has an extension in a complex domain denoted by the
same letter, and defined by
\begin{equation}
\label{saddle}
 df + \varepsilon \omega_\varepsilon = 0
\end{equation}
where $ \omega_\varepsilon  =P(x,y,\varepsilon) dx + Q(x,y,\varepsilon) dy$ is
a one-form, and $P,Q$ are real-analytic in $x,y,\varepsilon$ in a neighborhood
of $\Gamma_k$.

Parameterize the segment $\sigma$ by the "synchronized" local variable $t=
f|_\sigma$ and let $\gamma(t)\subset F_t, t>0$ be the continuous family of
periodic orbits of $X_0$ which tend to the polycycle $\Gamma_k\subset F_0$ as
$t$ tends to $0$. To the family $\{\gamma(t)\}$ we associate the trivial first
return map
$$
P_0: \sigma \rightarrow \sigma, P_0=id
$$
which allows an analytic continuation for $\varepsilon\neq 0$ to a first return
map
\begin{equation}\label{pmap}
\begin{array}{rcl}
P_\varepsilon: \sigma & \rightarrow &\sigma \\
     t &\mapsto & P_\varepsilon(t)= t + \varepsilon^d M_d(t) + \dots
\end{array}
\end{equation}
The dots above mean a function in $t,\varepsilon$ which, for every fixed $t$
such that $P_\varepsilon(t)$, $\varepsilon\sim 0$ is defined, is of the type
$0(\varepsilon^{d+1})$. The so called Poincar\'{e}-Pontryagin function function
$M_d$ may be explicitly computed, see \cite{pont34,fran96,fran98}.

More generally, let $\gamma(t) \subset F_t$ be any continuous family of closed
loops intersecting the cross-section $\sigma$. For $|\varepsilon|$ sufficiently
small we define in a similar way the  holonomy map
$$
h_{\gamma}^\varepsilon : \sigma \rightarrow \sigma
$$
related to the family of loops $\{ \gamma(t)\}_t$ and the deformed foliation
$\mathcal{F}_\varepsilon$.  By analogy to the Poincar\'{e} return map we have
$$
h_{\gamma}^\varepsilon (t)= t + \varepsilon^d M_d(t) + \dots
$$
\begin{equation}\label{pmap1}
\begin{array}{rcl}
h_{\gamma}^\varepsilon : \sigma & \rightarrow &\sigma \\
     t &\mapsto & h_{\gamma}^\varepsilon (t)= t + \varepsilon^d M_d(t) + \dots
\end{array}
\end{equation}
where the meaning of the dots is as before, and the number $d$ depends on
$\{\gamma(t)\}_t$ and $\mathcal{F}_\varepsilon$ .

The holonomy map $h_{\gamma}^\varepsilon$ depends on the choice of $\sigma$. In
contrast to this, the Poincar\'{e}-Pontryagin function $M_d$  does not depend on
the cross-section $\sigma$, it depends on the free homotopy class of the loop
$\gamma(t)$ only. Further, it can
 be
expressed in terms of iterated path integrals of length at most $d$, along
suitable meromorphic differential one-forms. It satisfies therefore a linear
differential equation which has a Fuchs type singularity at $t=0$, see
\cite{gavr05, gail05}. Thus, the leading term  of $M_d$ has the form
$$
t^p (\log t)^q .
$$
where $p$ is an eigenvalue of the indicial equation of the Fuchsian equation
related to the regular singular point $t=0$.
\begin{definition}
We shall call $p$ the \emph{characteristic number} of the holonomy map
$h_{\gamma}^\varepsilon$ and denote
$$
\nu(h_{\gamma}^\varepsilon)=p.
$$
\end{definition}
In the Hamiltonian case the number $\nu(P_\varepsilon)$ is rational, because
the corresponding monodromy operator is quasi-unipotent \cite{gano08}.

To formulate the main result of the paper consider, more specifically, the case
$k=2$ ( a double heteroclinic loop). As it follows from \cite{gail05}, the
function $M_d$ is in fact an Abelian integral and can be written as
$$
M_d(t)=\frac{1}{t^{d-1}}\int_{\gamma_0(t)} \tilde{\omega}
$$
for suitable analytic one-form $\tilde{\omega}$.

Let $\delta_1(t), \delta_2(t)$ be two continuous families of closed loops
vanishing at the saddle points $p_1$ and $p_2$. We suppose that orientations of
the loops "agree" in the sense that the intersection indices of the homotopy
classes of $\delta_1(t), \delta_2(t)$ with the homotopy class of the periodic
orbit $\gamma(t)$ is one and the same.

The cyclicity $Cycl(\Gamma_2,\mathcal{F}_\varepsilon)$ of the $2$-saddle cycle
$\Gamma_2$ with respect to the deformed foliation  $\mathcal{F}_\varepsilon$ is
the maximal number of limit cycles which bifurcate from $\Gamma_2$ near
$\varepsilon=0$, see~\cite{rous98} for a precise definition. An upper bound for
the cyclicity $Cycl(\Gamma_2,\mathcal{F}_\varepsilon)$ is given in terms of the
characteristic numbers of the holonomies associated to $\Gamma_2$ as follows
\begin{theorem}
\label{main}
\begin{equation}\label{mainbound}
Cycl(\Gamma_2,\mathcal{F}_\varepsilon) \leq 1+ \nu(P_\varepsilon)+
\max\{\nu(h^\varepsilon_{\delta_1}),\nu(h^\varepsilon_{\delta_2})\}+
\nu(h^\varepsilon_{\delta_1}\circ h^\varepsilon_{\delta_2}) .
\end{equation}
\end{theorem}
It is tempting to conjecture that in general the cyclicity
$Cycl(\Gamma_2,\mathcal{F}_\varepsilon)$ is bounded by a similar expression in
terms of the characteristic numbers of the holonomies associated to $\Gamma_k$.
Indeed, in the homoclinic case, $k=1$, by repeating the proof of
Theorem\ref{main} one obtains
$$
Cycl(\Gamma_1,\mathcal{F}_\varepsilon) \leq \nu(P_\varepsilon)+
\nu(h_{\delta_1}) .
$$
We have, typically
$$
P_\varepsilon(t)= t + \varepsilon M_1(t)+\dots
$$
and if $$ M_1(t) = f_1(t) \log(t) + f_2(t) \not\equiv 0, f_1(t)= O(t^{p}),
f_2(t)= O(t^{q})
$$
for some analytic functions $f_1, f_2$, then
\begin{equation*}
Cycl(\Gamma_1, \mathcal{F}_\varepsilon) \leq \min \{p,q \}+ q .
\end{equation*}
By the Roussarie's theorem \cite[Theorem C]{rous86} the exact upper bound in a
real domain in this case is $2p$ if $p<q$, and $2q-1$ if $p\geq q$. This
suggests that the bound of Theorem \ref{main} can be improved. In fact, the
bound (\ref{mainbound}) holds true for the number of complex limit cycles
accumulating on $\Gamma_2$ in a suitable neighborhood of it.
\\
\textbf{Example.} Suppose that $d=1$ in (\ref{pmap}), that is to say $M_1(t)=
\int_{\gamma(t)} \omega_0$ where  $\{\gamma(t)\}_{t>0}$ is the family of real
periodic orbits of $\mathcal{F}_0 = \{ df = 0\}$. Then we have
$$
M_1(t) = (f_1(t)   + f_2(t)) \log(t) + f_3(t)
$$
where the functions $f_1,f_2,f_3$ are analytic in a neighborhood of $t=0$,
$$
f_1(t) = \int_{\delta_1(t)} \omega_0, f_2(t) = \int_{\delta_2(t)} \omega_0
$$
$$
h_{\delta_1}(t)= t + \varepsilon f_1(t) +\dots, h_{\delta_2}(t)= t +
\varepsilon f_2(t) +\dots$$
$$
 h_{\delta_1}\circ h_{\delta_2} (t)= t +
\varepsilon (f_1(t)+f_2(t)) +\dots
$$
Suppose further that
$$
f_3(t)= O(t^p), f_1(t)= O(t^{p_1}), f_2(t)= O(t^{p_2}), f_1(t)+f_2(t) = O(t^q).
$$
Theorem \ref{main} implies that the cyclicity of $\Gamma_2$ is bounded by
$$
1+ \min \{p,q\} + \max \{p_1,p_2\} + q .
$$
In the case $p_1=p_2=q=p$ for instance, this gives
$$
Cycl(\Gamma_2,\mathcal{F}_\varepsilon) \leq 1+ 3 p .
$$
In this situation, and under the strong hypothesis that one of the connections
of $\Gamma_2$ remains unbroken, it has been proved in \cite[Theorem 8] {duro06}
that
$$
Cycl(\Gamma_2,\mathcal{F}_\varepsilon) \leq  2p-1 + \frac{p(p-1)}{2} .
$$

\section{The Dulac map}
\label{dulacsection}
\begin{figure}
\begin{center}
\resizebox{10cm}{!}{
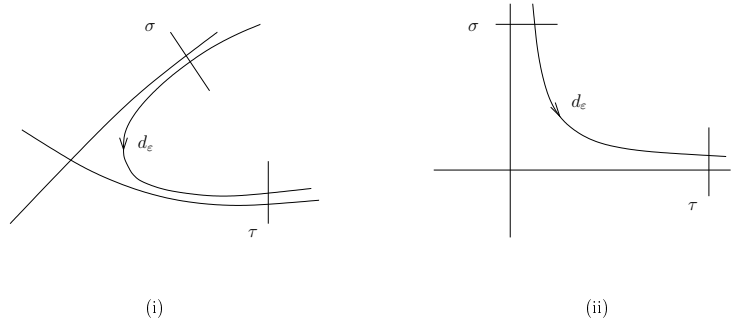}
\end{center} \caption{The Dulac map}\label{fig1a}
\end{figure}
Let $\mathcal{F}_\varepsilon$ be a real analytic foliation defined by
(\ref{saddle}) in a neighborhood of a hyperbolic Morse critical point of the
function $f$. For all sufficiently small $|\varepsilon|$ the foliation
$\mathcal{F}_\varepsilon$ has a singular point of saddle type, to which we
associate a Dulac map (or the transition map), as on fig.\ref{fig1a} (i). More precisely, for
all sufficiently small $|\varepsilon|$ the foliation $\mathcal{F}_\varepsilon$
has two separatrix  solutions, which are transversal analytic curves, depending
analytically on $\varepsilon$. We may suppose that they are the axes $\{x=0\}$
and  $\{y=0\}$ as on fig.\ref{fig1a} (ii). Let $\sigma, \tau$ be two complex cross-sections (complex
discs) to the two separatrices, parameterized by $z=f|_\sigma$ and $
z=f|_\tau$. In these coordinates the Dulac map is the germ of analytic map
$$
d_\varepsilon : (\mathbb{R}^+_*,0)\rightarrow (\mathbb{R}^+_*,0)
$$
defined as follows: if $z\in\sigma\cap \mathbb{R^+_*}$ then
$d_\varepsilon(z)\in \tau\cap \mathbb{R^+_*}$ is the intersection with $\tau$
of the orbit $\gamma_\varepsilon(z)$ of (\ref{saddle}), passing through $z\in
\sigma$. This geometric definition of $d_\varepsilon$ allows to control to a
certain extent its analytic continuation in a complex domain.
\subsection{Analytic continuation}
\label{scontinuation}
 The Dulac map is analytic and hence allows an analytic
continuation on some open subset of the universal covering  $\sigma_\bullet$ of
$\sigma\setminus \{0\}$. The domain of the  continuation depends on
$\varepsilon$, and obviously $d_0(z)\equiv z$.

Let us parameterize  the universal covering  $\sigma_\bullet$ by polar
coordinates $\rho > 0, \varphi \in \R$, $z= \rho \exp{\varphi}$.

\begin{theorem}
\label{continuation} There exists $\varepsilon_0>0$ and a continuous function
\begin{equation*}
\begin{array}{rcl}
\rho: \R& \rightarrow &\R^+_*\\
     \varphi &\mapsto & \rho(\varphi)
\end{array}
\end{equation*}
such that the Dulac map allows an analytic continuation in the domain
\begin{equation}\label{domain}
\{\varepsilon, \rho, \varphi)\in \C\times\sigma_\bullet: |\varepsilon|<
\varepsilon_0, 0< \rho < \rho(\varphi)\}
\end{equation}
\end{theorem}
The proof of the above Theorem in the $0$-parameter case is well known, and in
the multi-parameter case it is the same. For convenience of the reader it will
be given in Appendix \ref{A}. This proof shows even more: the analytic
continuation of the Dulac map in the domain (\ref{domain}) can be accomplished
in a geometric way as follows.

Let $\{\gamma_0(z)\}_z$, $\gamma_0(z)\subset F_z$ be a continuous family of
loops connecting $\sigma_\bullet$ to $\tau_\bullet$. For $z\in\sigma\cap
\mathbb{R^+_*}$ we suppose that $\gamma_0(z)$ is the real orbit of $df=0$
contained in the first quadrant $x\geq 0, y\geq 0$. We note that, although the
family $\{\gamma_0(z)\}_z$ is not unique, the relative homotopy class of each
loop $\gamma_0(z)$ is uniquely defined for all $z\in \sigma_\bullet$. It
follows from the proof of Theorem \ref{continuation} that $\{\gamma_0(z)\}_z$
allows a deformation to a family of paths $\{\gamma_\varepsilon(z)\}_z$,
 connecting $\sigma_\bullet$ to $\tau_\bullet$, tangent
to the leaves of $\mathcal{F}_\varepsilon$, and defined for all $\varepsilon,
z$ in the domain (\ref{domain}).

\subsection{The Poincar\'{e}-Pontryagin integral}
In what follows a crucial role will be played by the integral
$\int_{\gamma_0}\omega_0$, and its generalizations. Namely, let
$$
K\subset \{ (\rho, \varphi)\in \sigma_\bullet: 0< \rho < \rho(\varphi)\}
$$
be a compact set, where  $\rho(\varphi)$ is as in Theorem \ref{continuation}.
As $\int_{\gamma_{\varepsilon}(z)}\omega_\varepsilon$ is continuous in $z,
\varepsilon$ and
$$
d_\varepsilon(z)-z =
  \int_{\gamma_{\varepsilon}(z)} df = \varepsilon \int_{\gamma_{\varepsilon}(z)}
  \omega_\varepsilon
$$
then the following Lemma holds
\begin{lemma}[Pontryagin\cite{pont34}]
\begin{equation}\label{pontryagin1}
d_\varepsilon(z)=z + \varepsilon \int_{\gamma_{0}(z)}
  \omega_0 + O(\varepsilon^2)
\end{equation}
uniformly in $z\in K$.
\end{lemma}
The function $\int_{\gamma_{0}(z)} \omega_0$ is the so called
Poincar\'{e}-Pontryagin integral associated to the
  deformed foliation $\mathcal{F}_\varepsilon$. It follows from the argument
  principle that, that if $|\varepsilon|$ is sufficiently small,
   the number of the zeros  of $d_\varepsilon(z)-z$  in the compact $K$ is bounded by the
   number of the zeros of the Poincar\'{e}-Pontryagin integral $\int_{\gamma_{0}(z)}
  \omega_0$ in $K$ (counted with multiplicity).
It might happen, however, that the Poincar\'{e}-Pontryagin integral vanishes
identically. In all cases there is an integer $d\geq 1$ and an analytic
function $M_d \neq 0$ in a neighborhood of $K$, such that
\begin{equation}\label{pontryagin2}
d_\varepsilon(z)=z + \varepsilon^d M_d(z) + O(\varepsilon^{d+1})
\end{equation}
uniformly in $z\in K$, provided that the Dulac map is not the identity map.
$M_d$ is the so called higher order Poincar\'{e}-Pontryagin function and its zeros
in $K$ bound as before the number of the zeros of $d_\varepsilon(z)-z$. As we
already mentioned in section \ref{description}, there is an integral
representation for $M_d$ as an iterated integral of length at most $d$ along
$\gamma_0(z)$.

Our aim  is to obtain a bound for the zeros of $d_\varepsilon(z)-z$ in a domain
$K$ which is open and connected. Even if the estimate (\ref{pontryagin1}),
(\ref{pontryagin2}) do allow an extension to such a domain $K$, the argument
principle can not be directly used. For this purpose we consider rather the
imaginary part of the Dulac map.

\subsection{The zero locus of the imaginary part of the Dulac map}
\label{zerolocus}
 We shall describe the zero locus of the imaginary part of the Dulac map
$d_\varepsilon$ in an appropriate sector
\begin{equation}\label{zero}
\mathcal{H}_\varepsilon = \{z\in \C: \Im d_\varepsilon(z) = 0\} \cap
\mathcal{D}
\end{equation}
$$
\mathcal{D}= \{ 0<\rho<\rho(\varphi), 0<\varphi< \frac{3\pi}{2} \} .
$$
The surprising fact about $\mathcal{H}_\varepsilon$ is that it is a smooth
real-analytic plane curve in $\mathcal{D} \subset \mathbb{R}^2=\mathbb{C}$.
Even better, the curve  $\mathcal{H}_\varepsilon$ can be conveniently
approximated in terms of higher order Poincar\'{e}-Pontryagin functions.

The foliation $\mathbb{F}_0$ has a first integral defining a fibration with
fibers $F_t$, see (\ref{fiber}). Let $\delta(t)\subset F_t$ be a continuous
family of closed loops $\delta(t)\subset F_t$ vanishing at the saddle point
when $t$ tends to $0$. The orientation of $\delta(t)$ is chosen as follows. Let
$\gamma_0(t)$ be the family of loops defined in the Appendix \ref{A}. For real
positive $t$ they coincide with real orbits of $\mathcal{F}_0$ connecting
$\sigma^+$ to $\tau^+$. Then, the homotopy classes of $\gamma_0, \delta$
satisfy
\begin{equation}\label{orientation} \gamma_0(te^{i\pi})- \gamma_0(te^{-i\pi}) =
\delta(t) .
\end{equation}

Therefore, the exact orientation of $\delta(t)$ can be computed by the
Picard-Lefschetz formula (but we do not need this).

Let $\tau$ be, as before, a cross section to the fiber $F_0$, see fig.
\ref{fig1a}. Consider the holonomy map $h_{\delta}^\varepsilon$ associated to
the family $\{\delta(t)\}_t$ and to the deformed foliation
$\mathcal{F}_\varepsilon$
\begin{equation}\label{vanishing}
\begin{array}{rcl}
h_{\delta}^\varepsilon : \tau & \rightarrow& \tau \\
     z&\mapsto&  z + \varepsilon^d M_d(z) + \dots
\end{array}
\end{equation}
 The anti-holomorphic involution
$$
(x,y)\mapsto (\bar{x},\bar{y})
$$
induces, for $t\in\R$ an anti-holomorphic involution
$$
F_t \rightarrow F_t
$$
which on its turn sends the free homotopy class of the loop $\delta(t)$ to the
class of $-\delta(t)$. Therefore the function $M_d$ is pure imaginary for real
values of $t$.

\begin{lemma}
\label{mainlemma}
The zero locus $\mathcal{H}_\varepsilon$ of the imaginary
part of the Dulac map is a smooth real-analytic curve of $\R^2=\C$ of the form
\begin{equation}\label{sub}
\mathcal{H}_\varepsilon = \{z=u+iv: v= \frac{\varepsilon^d}{2i} M_d(u) +
\varepsilon^{d+1} R(u,\varepsilon), u< 0\} \cap \mathcal{D}
\end{equation}
where $R(u,\varepsilon)$ is an analytic function.
\end{lemma}
\begin{figure}
\begin{center}
\resizebox{7cm}{!}{
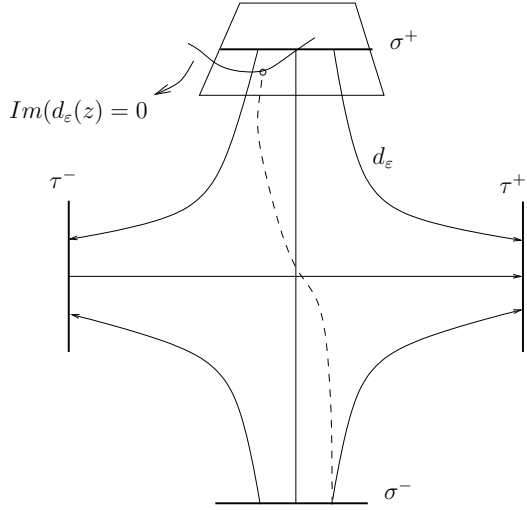}
\end{center} \caption{The zero locus of the imaginary part of the Dulac map}
\label{dulac}
\end{figure}
The above Lemma is the main technical result of the present paper. The
analyticity of the zero locus $\mathcal{H}_\varepsilon$ is responsible for the
algebraic-like behavior of the Dulac map. \\
\textbf{Proof of Lemma \ref{mainlemma}}. Consider the cross-sections (complex
discs transversal to the separatrices) $\sigma^\pm$, $\tau^\pm$ as shown on
fig.\ref{dulac}, simultaneously parameterized as before by the restriction $z$
of the first integral $f(x,y)$ on them. The cross-sections $\sigma, \tau$ shown
on fig.\ref{fig1a} are  denoted, from now on, by $\sigma^+,\tau^+$. Denote
$$
\sigma^+_{\geq 0}= \sigma^+\cap \{ (x,y): f(x,y)\geq 0\}, \sigma^+_{\leq 0}=
\sigma^+\cap \{ (x,y): f(x,y)\leq 0\} \mbox{  etc.}
$$
Let   $\{\gamma_\varepsilon(z)\}_z$ be the continuous family of paths, defined
in Appendix \ref{A}. The point $z\in\sigma^+$ belongs to the zero locus
$\mathcal{H}_\varepsilon$ if and only if the end of the path
$\gamma_\varepsilon(z)$ belongs to $\tau^+_{<0}$. Therefore, such a path allows
a decomposition in a product
$$
\gamma_\varepsilon(z) = \alpha_\varepsilon(z) \circ \beta_\varepsilon(z)
$$
where $\alpha_\varepsilon(z)$ is a path connecting $z\in\sigma^+$ to a point on
$\sigma^-_{<0}$ and $\beta_\varepsilon(z)$ is a path connecting the latter
point to a point on $\tau^+_{<0}$, where $\beta_\varepsilon(z)\subset \R^2$. It
follows that $\mathcal{H}_\varepsilon$ is the image of $\sigma^-_{<0}$ under
the holonomy map
$$
h_{\alpha_0^{-1}}^\varepsilon : \sigma^-_{<0} \rightarrow \sigma ^+ .
$$
This already proves that the closure of $\mathcal{H}_\varepsilon\subset \R^2$
is a smooth analytic curve. Once having said this, it is clear that
$\mathcal{H}_\varepsilon$ can be conveniently parameterized, which we do next.

 To prove $(\ref{sub})$, let us note that for $t\in\R$
the two ends of the loop $\alpha_0(t)$ are real and hence this also holds true
for the complex-conjugate loop $\overline{\alpha_0(t)}$. The loop
$\alpha_0(t)\circ \overline{\alpha_0^{-1}(t)}$ is therefore closed and is
homotopic to $\delta(t)\subset F_t$ defined above.

More explicitely
\begin{eqnarray*}
 2 i v & = & h_{\alpha_0^{-1}}^\varepsilon(t) -
\overline{h_{\alpha_0^{-1}}^\varepsilon(t)}\\
   &=& h_{\alpha_0^{-1}}^\varepsilon(t) -
h_{\overline{\alpha_0^{-1}}}^\varepsilon(t)\\
& = &( h_{\alpha_0^{-1}}^\varepsilon \circ h_{\overline{\alpha_0}}^\varepsilon
- id) \circ h_{\overline{\alpha_0^{-1}}}^\varepsilon(t)\\
& = &( h_{\alpha_0^{-1}\circ \overline{\alpha_0}}^\varepsilon
- id) \circ h_{\overline{\alpha_0^{-1}}}^\varepsilon(t)\\
& = & (h_\delta^\varepsilon - id)\circ
h_{\overline{\alpha_0^{-1}}}^\varepsilon(t) .
\end{eqnarray*}
This, together with (\ref{vanishing}) and
$$
h_{\overline{\alpha_0^{-1}}}^\varepsilon(t) = t + O(\varepsilon),
h_{\alpha_0^{-1}}^\varepsilon(t) = t + O(\varepsilon)
$$
implies
\begin{eqnarray*}
  u &=& t + O(\varepsilon) \\
  v &=& \frac{\varepsilon^d}{2i}( M_d(t) + O(\varepsilon))
\end{eqnarray*}
where, by abuse of notation,  $O(\varepsilon)$ means a function analytic in
$t,\varepsilon$, which vanishes identically for $\varepsilon = 0$. This proves
the identity (\ref{zero}).$\triangle$

\section{Cyclicity of two-saddle cycles}
\label{mainsection}
\begin{figure}
\begin{center}
\resizebox{7cm}{!}{
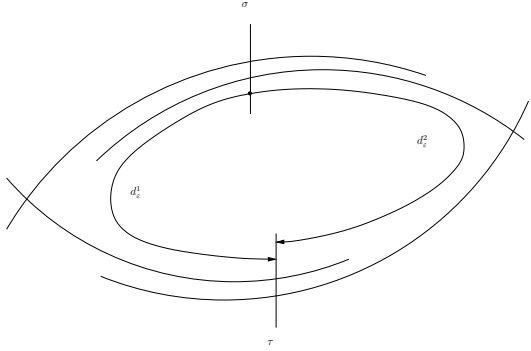}
\end{center} \caption{The Dulac maps $d^1_\varepsilon$ and
$d^2_\varepsilon$}\label{fig4}
\end{figure}
\begin{figure}
\begin{center}
\resizebox{7cm}{!}{
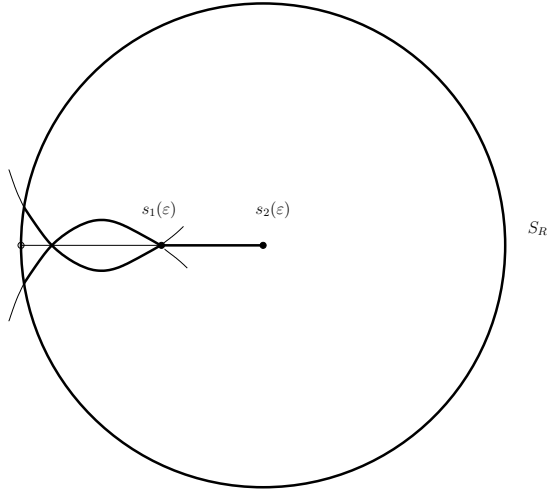}
\end{center} \caption{The domain $\mathcal{D}_\varepsilon$}\label{fig5}
\end{figure}

In this section we prove Theorem \ref{main}. Using the notations introduced in
section \ref{description},  we suppose that the vector field $X_0$ has a
two-saddle loop $\Gamma_2$ and an analytic first integral   $f$ in a
neighborhood of it. We suppose that $f$ has Morse critical points at the  two
saddle points $p_1$, $p_2$ of $X_0$. Consider the Dulac maps $d^1_\varepsilon$,
$d^2_\varepsilon$ associated to the corresponding foliation, see
fig.\ref{fig4}. Each map $d^i_\varepsilon$ is a composition of a "local" Dulac
map (as in section \ref{dulacsection}) and two holomorphic holonomy maps. From
this it follows that Lemma \ref{mainlemma} applies to  $d^i_\varepsilon$,
$i=1,2$, too. We parameterize each cross-section by the restriction $z$ of $f$
on it. The function $d^i_\varepsilon$, $i=1,2$, is multivalued and has a
critical point at $s_i(\varepsilon)\in \mathbb{R}$, $s_i(0)=0$. The functions
$s_i$ are real analytic.  We consider first the case $\varepsilon>0$ and we may
suppose that $s_i(0)=0$, $s_1(\varepsilon)< s_2(\varepsilon)$ for all
sufficiently small $\varepsilon$, (the case $\varepsilon < 0$ is studied in the
same way). Our aim is to bound the number of those zeros of the displacement
map $ d^1_\varepsilon - d^2_\varepsilon $ which are real, bigger than
$s_2(\varepsilon)$ and tend to $0$ as $\varepsilon$ tends to $0$. Note that
these zeros correspond to the fixed points of the Poincar\'{e} first return map
$P_\varepsilon$. Indeed,
\begin{equation}\label{dmap}
d^1_\varepsilon - d^2_\varepsilon = d^2_\varepsilon \circ ( (P_\varepsilon - id
),\mbox{  where  } d^2_0= id .
\end{equation}
We shall count the zeros of the displacement map in the larger complex domain
 $\mathcal{D}_\varepsilon$ of the universal covering
of $\mathbb{C}\setminus \{s_1(\varepsilon)\}$ defined as follows. It is bounded
by the circle
\begin{equation}\label{sr}
S_R= \{z: |z| = R \},
\end{equation}
by the interval $[s_1(\varepsilon),s_2(\varepsilon)]$, and by the zero locus of
the imaginary part of the Dulac map $d^1_\varepsilon$ for $\Re (z) <
s_1(\varepsilon)$, as it is shown on fig.\ref{fig5}. The numbers $\varepsilon,
R$ are subject to certain conditions explained bellow. The zeros of an analytic
function in a complex domain equal the increase of the argument of the function
along the border of the domain, divided by $2\pi$ (the argument principle). To
bound the increase of the argument we shall count the number of the zeros of
the imaginary part of the function, along the border of the domain.

Choose first the real numbers $\varepsilon_0, R >0$ as follows. Let
$$
d^1_\varepsilon(z) - d^2_\varepsilon(z) = \varepsilon^d M_d(z)
+0(\varepsilon^{d+1})
$$
and let $z^\nu (\log z)^\mu$ be the leading term of $M_d(z)$. Then, by
(\ref{dmap}), $\nu=\nu(P_\varepsilon)$ is the characteristic number of the
Poincar\'{e} map $P_\varepsilon$ associated to $\Gamma_2$. We choose $R>0$ so
small, that the increase of the argument of $M_d(z)$ along the circle $S_R$ is
sufficiently close to the increase of the argument of $z^\nu (\log z)\mu$ along
$S_R$. We fix $R$ and choose $\varepsilon_0>0$ so small with respect to $R$,
that for all $\varepsilon$, $|\varepsilon|<\varepsilon_0$, the increase of the
argument of $d^1_\varepsilon(z) - d^2_\varepsilon(z)$ along the circle $S_R$ is
sufficiently close to the increase of the argument of $M_d(z)$ along $S_R$.
This is indeed possible, according to Lemma \ref{pontryagin1}. The conditions
that we impose on $\varepsilon_0,\varepsilon$ and $ R$ will be denoted (by
abuse of notations) as follows
$$
1 >> R >> \varepsilon_0 > \varepsilon  >0 .
$$

To evaluate the increase of the argument of $ d^1_\varepsilon(z) -
d^2_\varepsilon(z)$ along the interval $[s_1(\varepsilon),s_2(\varepsilon)]$,
we bound the zeros of its imaginary part which equals  (along the interval
$[s_1(\varepsilon),s_2(\varepsilon)]$) to the imaginary part of
$-d^2_\varepsilon(z)$. In other words, we need to estimate the number of
intersection points (counted with multiplicity) between the zero locus of the
imaginary part of $d^2_\varepsilon(z)$ and
$[s_1(\varepsilon),s_2(\varepsilon)]$. According to Lemma \ref{mainlemma} this
number of intersection points is bounded by the multiplicity of the zero at the
origin of the Poincar\'{e}-Pontryagin function of the holonomy map
$h^\varepsilon_{\delta_2}$. This multiplicity equals
$\nu(h^\varepsilon_{\delta_2})$.

Finally, we arrive at the most delicate point in the proof of Theorem
\ref{main} :  evaluate the increase of the argument of $ d^1_\varepsilon(z) -
d^2_\varepsilon(z)$ along the zero locus of the imaginary part of the Dulac map
$d^1_\varepsilon$ for $\Re (z) < s_1(\varepsilon)$. For this purpose we bound
the zeros of the imaginary part  $ \Im{(d^1_\varepsilon(z) -
d^2_\varepsilon(z))}$, along the zero locus of $\Im d^2_\varepsilon(z)$. Thus,
we need to estimate the number of intersection points (counted with
multiplicity) between the zero locus of the imaginary part of
$d^2_\varepsilon(z)$, and the zero locus of the imaginary part of
$d^1_\varepsilon(z)$.

Recall that to a Dulac map $d^i_\varepsilon$ we associated a family of
vanishing loops $\delta_i(z)$ with orientation prescribed by
(\ref{orientation}). With this convention, the orientation of the loops
$\delta_1$ and $\delta_2$ do not agree : if $\gamma(t)$ is the family of
periodic orbits of $\mathcal{F}_0$ for $t>0$, then the intersection indices of
the homotopy classes of $\delta_1(t), \delta_2(t)$ with $\gamma(t)$ have
opposite signs. In order to have the same convention as in the formulation of
Theorem \ref{main} we reverse the orientation of $\delta_2$. With this
convention if
\begin{equation}\label{vanishing2}
\begin{array}{rcl}
h_{\delta_i}^\varepsilon :
     z&\mapsto&  z + \varepsilon^{d^{i}} M_d^{i}(z) + \dots
\end{array}
\end{equation}
then, by Lemma \ref{mainlemma}, the zero locus of the imaginary part of the
holonomy maps $h_{\delta_1}^\varepsilon, h_{\delta_2}^\varepsilon$ is given by
\begin{equation}\label{sub1}
 \{z=u+iv: v= \frac{\varepsilon^{d^1}}{2i} M_{d^{1}}(u) +
\varepsilon^{d^{1}+1} R^1(u,\varepsilon), u< 0\} \cap \mathcal{D}
\end{equation}
\begin{equation}\label{sub2}
 \{z=u+iv: v= -\frac{\varepsilon^{d^{2}}}{2i} M_{d^{2}}(u) +
\varepsilon^{d^{2}+1} R^2(u,\varepsilon), u< 0\} \cap \mathcal{D}
\end{equation}
respectively, where $R^1, R^2$ are appropriate analytic functions.

 We conclude that the number of intersection points of the above analytic curves coincides with the
 multiplicity of the zero at the origin of either $M_{d^{1}}$ (if $d^{1}< d^2$),
 or $M_{d^{2}}$ (if $d^2<d^1$), or $M_{d^{1}} + M_{d^{2}}$ (if $d^1=d^2$). The number of intersection
 points equals therefore to the characteristic number $\nu(h^\varepsilon_{\delta_1}\circ h^\varepsilon_{\delta_2})$.

Summing up the above information we get that the increase of the argument of $
d^1_\varepsilon(z) - d^2_\varepsilon(z)$ along the boundary of the complex
domain $\mathcal{D}_\varepsilon$ is not bigger than
$$
\nu(P_\varepsilon) + \nu(h^\varepsilon_{\delta_2}) +
\nu(h^\varepsilon_{\delta_1}\circ h^\varepsilon_{\delta_2}) +2 .
$$
The above estimate can be slightly improved, by taking into consideration the
fact that the imaginary part of $d^1_\varepsilon(z)-d^2_\varepsilon(z)$
vanishes at $s_2(\varepsilon)$. Theorem \ref{main} is proved.$\Box$
\def\cprime{$'$} \def\cprime{$'$} \def\cprime{$'$} \def\cprime{$'$}
  \def\cprime{$'$}
\appendix
\section{Proof of Theorem \ref{continuation}}
\label{A}
\begin{figure}
\begin{center}
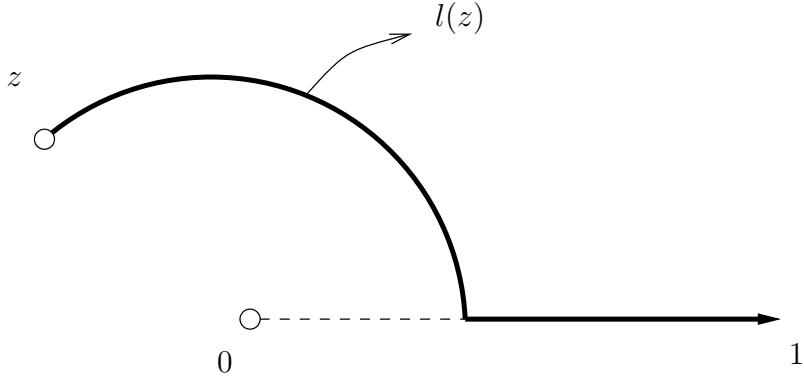
\end{center} \caption{The path $l(z)$}\label{fig6}
\end{figure}
In this appendix we will prove Theorem \ref{continuation} in the slightly more
general context of multi-parameter analytic deformations. This will be used in
Appendix \ref{B}. Note that the zero-parameter case is well known
\cite{lora05,rous98}.

 Consider a $N$-parameter analytic family of analytic plane vector
fields $X_\lambda$, $\lambda \in(\mathbb{R}^N,0)$, such that $X_0$ has a
hyperbolic singular point at the origin. It is known, since Briot and Bouquet,
that $X_0$ has two transversal invariant analytic curves which can be supposed
to coincide with the axes $x=0$ and $y=0$, that is to say
$$
X_0= \lambda_1 x(1+...) \frac{\partial}{\partial x} +  \lambda_2 y(1+
...)\frac{\partial}{\partial y}, \; \lambda_1\lambda_2 < 0 .
$$
The proof is as follows : a formal change of the variables first removes some
(but not all) non-resonant terms, and then one verifies the convergency of the
transformation, see \cite[Appendice II]{mamo80}. Exactly the same proof
applies, however, to the family $X_\lambda$. One can show  in this way that
$X_\lambda$ is analytically orbitally equivalent to the following (slightly
improved) normal form
\begin{equation}\label{a1}
x\frac{\partial}{\partial x} +  y(r+ xy. a)\frac{\partial}{\partial y}
\end{equation}
where $r=r(\lambda)$, $a=a(x,y,\lambda)$ are appropriate analytic functions in
their arguments, $r(0)<0$, see \cite[Appendice 1]{rous89}. We shall suppose,
without loss of generality, that there exists a constant $c>0$ such that $r, a$
are analytic in the complex domain
$$
\mathcal{D}_c= \{x,y,\lambda): |x|<2, |y|< 2, |\lambda| < c \}.
$$
After a further linear re-scaling of $x,y$, we may suppose that
$|a(x,y,\lambda)|$ is so small in $\mathcal{D}_c$, that
$$
r(\lambda)+ xy. a(x,y,\lambda) \neq 0 \; .
$$
After this preparation, choose the cross-sections
 $$\sigma=\{y=1\},
\tau=\{x=1\} $$
 and consider the corresponding Dulac map
\begin{equation*}
\begin{array}{rcl}
d_\lambda : \sigma    & \rightarrow & \tau \\
     z &\mapsto & d_\lambda(z) \; .
\end{array}
\end{equation*}
To prove Theorem \ref{continuation} we have to show that the constant $c>0$ can
be chosen in such a way, that for every $\varphi_0>0$ there exists $0<z_0<1$,
such that the Dulac map allows an analytic continuation in the sector
$$
\{ z\in \mathbb{C}: |z|< z_0, |arg(z)|< \varphi_0\}\times \{\lambda : |\lambda|
< c\} .
$$
The proof  is similar to the proof of \cite[Theorem 7]{rous98}, the only
difference being the presence of the parameter $\lambda$. We shall construct a
continuous family of paths $\gamma_\lambda(z)$ contained in the leaves of the
foliation $\mathcal{F}_\lambda$ defined by the vector field $X_\lambda$ in
$\mathbb{C}^2$. Each path $\gamma_\lambda(z)$ starts at the point $(x=z,y=1)$
and ends at the point $(x=1, y=d_\lambda(z)$. The path $\gamma_\lambda(z)$ is
constructed by lifting the path $l(z)$ contained in the $x$-plane $\{y=0\}$ and
shown on fig.\ref{fig6}, with respect to the projection
\begin{equation*}
\begin{array}{rcl}
\pi  : \mathbb{C }^2 \;\;\;  & \rightarrow & \mathbb{C } \\
     (x,y) &\mapsto & x
\end{array}
\end{equation*}

\begin{figure}
\begin{center}
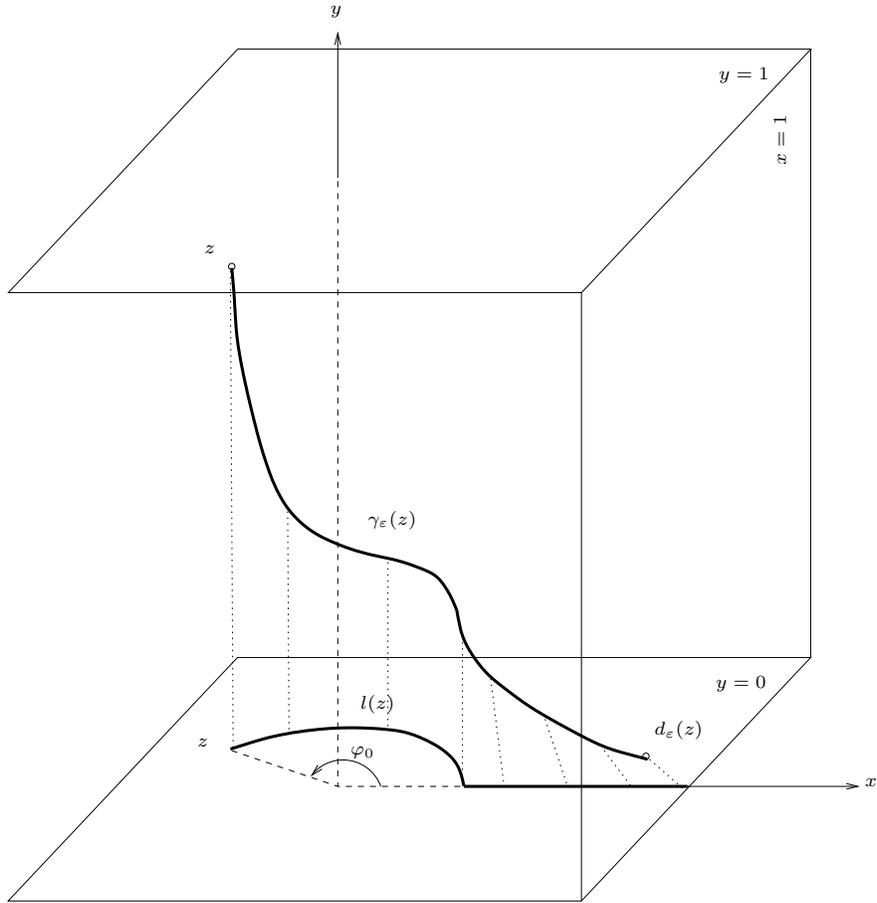
\end{center}
\caption{The path $\gamma_\lambda(z)$} \label{gz}
\end{figure}

 Indeed, the foliation $\mathcal{F}_\lambda$ is transversal to the
projection $\pi$ except along the leaf $x=0$, provided that $r+ xy. a \neq 0$.
The resulting path path $\gamma_\lambda(z)$ is shown on fig. \ref{gz}.

 To prove the existence of $\gamma_\lambda(z)$, consider the  solution $y=y(x)$
 associated to the vector field $X_\lambda$,
with initial condition $y(z)=1$. We have to show that the solution $y=y(x)$
exists when $x$ is restricted to the path $l(z)$. The path $l(z)$ is composed
by an arc and a segment. We consider them separately
\begin{itemize}
    \item
    Along the arc
$$
x= |z| e^{i\varphi}, y= 0, 0<\varphi < arg(z)
$$
parameterized by $\varphi$ we have
\begin{eqnarray*}
  dy &= &- y (r + |z| e^{i\varphi} y .a) d\varphi\\
 d|y| &= &- |y|.|z| \Im (e^{i\varphi} y .a) d\varphi .
\end{eqnarray*}
Therefore, if $|z|<z_0$ is sufficiently small, then $|y(|z| e^{i\varphi})|<2$
when $0<\varphi < arg(z)$
    \item
Along the segment
$$
x\in [|z|,1]
$$
we have similarly
\begin{eqnarray*}
 x dy &= &- y (r + x y .a) dx\\
 x d|y| &= &- |y|(r + \Re (x y .a)) dx .
\end{eqnarray*}
The derivative $\frac{d|y|}{dx} $ is therefore negative, the function $|y|(x)$
decreasing, so $|y|(x) < 1$.
\end{itemize}
$\Box$

\section{Multi-parameter deformations of Hamiltonian two-saddle loops}
\label{B}
 Consider, as in the preceding Appendix, a $N$-parameter analytic
family of analytic plane vector fields $X_\lambda$, $\lambda
\in(\mathbb{R}^N,0)$. We suppose that $X_0$ has a Hamiltonian two-saddle loop
$\Gamma_2$ bounding a period annulus. This case is easily reduced to the
one-parameter case studied in the present paper by making use of a standard
procedure based on the Hironaka desingularization theorem. In this Appendix we
indicate the main steps.

\subsection{Principalization of the Bautin ideal}
\label{b1}
 In this section we follow \cite{gavr08, rous01}. Let $z_0\in
\sigma$, $z_0\not\in \Gamma_2$ and consider the Poincar\'{e} map
$$
P_\lambda(z)= z + \sum_{i=1}^\infty a_i(\lambda)(z-z_0)^i .
$$
The Bautin ideal, associated to $P_\lambda$, is the ideal $\mathcal{I}=<a_i>$,
generated by the germs of the analytic coefficients $a_i(.)$. It is Noetherian,
so generated by a finite number of coefficients, and moreover does not depend
on the choice of $z_0$ \cite{rous98}. More generally, let $\{\gamma(z)\} $ be
any continuous family of closed loops in the fibers of the foliation
$\mathcal{F}_0$, intersecting the cross-section $\sigma$. For all sufficiently
small $\lambda$ the (germ of) holonomy map
$$
h^\lambda_\gamma : \sigma \rightarrow \sigma
$$
is defined. In the same way, we associate to $h^\lambda_\gamma(z)$ a Bautin
ideal, generated by the coefficients $a_i(\lambda)$ of the expansion of
$h^\lambda_\gamma(z)$ with respect to $z-z_0$. As before it is Noetherian, and
does not depend on the choice of $z_0$ (with the same proof).

We may assume, without any loss of generality, that the Bautin ideal is
principal. For this we use a variant of Hironaka's desingularization theorem as
follows.

 Let $\varphi_0,\varphi_1,\dots,\varphi_p$ be non-zero analytic functions on a smooth
complex or real analytic variety $X$.
The indeterminacy points of the rational map
$$
\varphi: X \dashrightarrow \P^p
$$
can be eliminated as follows \cite{hiro64,bimi89}

\begin{theorem}[Hironaka desingularization]
\label{hiroth}
 There exists a smooth analytic variety $\tilde{X}$ and a proper analytic map $\pi:
\tilde{X}\rightarrow X$ such that the induced map $\tilde{\varphi}=\varphi\circ
\pi$ is analytic.
$$
 \xymatrix{
    \tilde{X}  \ar[d]_\pi \ar[rd]^{\tilde{\varphi}}&  \\
    X \ar@{-->}[r]^\varphi & \P^p }
  $$
  \end{theorem}
 Let $\mathcal{O}_X$ be the sheaf of analytic functions on $X$ and consider the ideal
sheaf $I \subset \mathcal{O}_X$ generated by
$\varphi_0,\varphi_1,\dots,\varphi_p$. The inverse image ideal sheaf of $I$
under the map $\pi: \tilde{X}\rightarrow X$ will be denoted $\pi^*I$. This is
the ideal sheaf generated by the pull-backs of local sections of $I$. We note
that $\pi^*I$  may differ from the usual sheaf-theoretic pull-back, also
commonly denoted by $\pi^*I$. A simple consequence of Theorem \ref{hiroth} is
the following

\begin{corollary}
\label{hirocor} The inverse image ideal sheaf $\pi^*I$ is principal.
 \end{corollary}
This is called the principalization of $I$. Indeed, as the induced map
$\tilde{\varphi}$ is analytic, then for every $\tilde{\lambda}\in \tilde{X}$
there exists  $j$, such that the functions $\tilde{\varphi}_i/
\tilde{\varphi}_j$, $i=1,2,\dots,p$, are analytic in a neighborhood of
$\tilde{\lambda}$. Therefore there is a neighborhood $\tilde{U}$ of
$\tilde{\lambda}$ such that $\tilde{\varphi}_j|_{\tilde{U}}$ divides
$\tilde{\varphi}_i|_{\tilde{U}}$ in the ring of sections
$\mathcal{O}_{\tilde{U}}$ of the sheaf $\mathcal{O}_{\tilde{X}}$, that is to
say $I_{\tilde{U}}$ is generated by $\tilde{\varphi}_j|_{\tilde{U}}$.

In our context $X=\sigma$ is the cross-section to the family of periodic orbits
$\{\gamma(z)\}$ and $\varphi_0,\varphi_1,\dots,\varphi_p$ are the germs of
analytic functions which generate the Bautin ideal associated to the holonomy
map $h^\lambda_\gamma$. To apply Theorem \ref{hiroth} we assume that $\sigma$
is a polydisc on which $\varphi_i$ are analytic, and the divisors $(\varphi_i)$
intersect transversally the boundary of $\sigma$. After applying the Hironaka's
theorem, the origin $0$ of $\sigma$ is replaced by a compact divisor
$\pi^{-1}(0)$, along which the inverse image ideal sheaf $\pi^*I$ is principal.

Suppose now that we have  a holonomy map $h^\lambda_\delta$ associated to
another family of periodic orbits $\delta(t)$. Let $\psi_0, \psi_1,\dots ,
\psi_q$ be generators of the corresponding Bautin ideal $J$. As before we
assume that $\psi_i$ are analytic on $\sigma$, with divisors transversal to its
boundary. Applying twice the Hironaka's theorem we get a new analytic variety
smooth analytic variety $\tilde{X}$ and a proper analytic map $\pi:
\tilde{X}\rightarrow X$ such that the induced maps
$\tilde{\varphi}=\varphi\circ \pi$ and $\tilde{\psi}=\psi\circ \pi$ are
analytic, see diagram (\ref{pih}).
\begin{equation}\label{pih}
\xymatrix{
  & \ar[ld]_{\tilde{\psi}}  \tilde{X}  \ar[d]_\pi \ar[rd]^{\tilde{\varphi}} &  \\
 \P^q &   \ar@{-->}[l]_\psi X \ar@{-->}[r]^\varphi & \P^p }
\end{equation}
The inverse image ideal sheaf $\pi^*I$ and $\pi^*J$ are both principal along
the compact divisor $\pi^{-1}(0)$.

\subsection{Multi-parameter version of Theorem \ref{main}}
\label{b2}
 We begin by formulating the  multi-parameter version of
Lemma~\ref{mainlemma}, let $\{\delta(t) \}$  be the family of vanishing loops
defined in section~\ref{zerolocus}. In agreement with the preceding section,
let us suppose that the Bautin ideal of the holonomy map $h^\lambda_\delta$ is
principal. We have, therefore (compare to (\ref{vanishing})
\begin{equation}\label{vanishing1}
\begin{array}{rcl}
h_{\delta}^\lambda : \tau & \rightarrow& \tau \\
     z&\mapsto&  z + \tilde{\varphi}(\lambda)( \tilde{M}(z) + \tilde{R}(z,\lambda))
\end{array}
\end{equation}
where $\tilde{R}(.,.)$ is analytic, $\tilde{R}(z,0)=0$, and $\tilde{M}(.)$ is
the highest order Poincar\'{e}-Pontryagin function.
\begin{lemma}
\label{mainlemma1} The zero locus $\mathcal{H}_\varepsilon$ of the imaginary
part of the Dulac map is a smooth real-analytic curve of $\R^2=\C$ of the form
\begin{equation}\label{sub1+}
\mathcal{H}_\lambda = \{z=u+iv: v= \frac{\tilde{\varphi}(\lambda)}{2i}
(\tilde{M}(u) + R(u,\lambda)), u< 0\} \cap \mathcal{D}
\end{equation}
where $R(u,\lambda)$ is an analytic function, $R(u,0)=0$.
\end{lemma}
The proof of the above Lemma is completely analogous to that of
Lemma~\ref{mainlemma} and is therefore omitted.

In  the proof of Theorem \ref{main} we used four Bautin ideals associated to
the holonomy maps
\begin{equation}\label{bautin4}
P_\lambda=h^\lambda_\gamma, h^\lambda_{\delta_1}, h^\lambda_{\delta_2},
h^\lambda_{\delta_1}\circ h^\lambda_{\delta_2}
\end{equation}
where $\{\gamma(z)\}$ is the family of periodic orbits associated to the
annulus, $\{\delta_1(z)\}, \{\delta_2(z)\}$ are the vanishing families of loops
associated to the saddle points. After an appropriate blow up $\pi$ we may
suppose that the inverse image ideal sheafs of the corresponding four Bautin
ideals are principal along the compact divisor $\pi^{-1}(0)\subset \tilde{X} $,
see section~\ref{b1}. Let $\tilde{\lambda}$ be a local variable on the smooth
variety $\tilde{X}$. The cyclicity $
Cycl(\Gamma_2,(\mathcal{F}_{\tilde{\lambda}}, \mathcal{F}_{\tilde{\lambda}_0}))
$ is the maximal number of limit cycles which $\mathcal{F}_{\tilde{\lambda}}$
can have in an arbitrarily small neighborhood of $\Gamma_2$, when
$\tilde{\lambda}$ tends to $\tilde{\lambda}_0$. Denote also
$$
Cycl(\Gamma_2,\mathcal{F}_\lambda) =
Cycl(\Gamma_2,(\mathcal{F}_\lambda,\mathcal{F}_0) .
$$
Clearly
$$
Cycl(\Gamma_2,\mathcal{F}_\lambda)= \sup_{\tilde{\lambda}\in\pi^{-1}(0)}
Cycl(\Gamma_2,(\mathcal{F}_{\tilde{\lambda}}, \mathcal{F}_{\tilde{\lambda}_0}))
$$
and because of the compactness of $\pi^{-1}(0)$, there exists
$\tilde{\lambda}_0 \in \pi^{-1}(0)$ such that
$$
Cycl(\Gamma_2,\mathcal{F}_\lambda) =
Cycl(\Gamma_2,(\mathcal{F}_{\tilde{\lambda}}, \mathcal{F}_{\tilde{\lambda}_0}))
.
$$

The above considerations show that, without any harm, we may suppose that
$\tilde{X}=\sigma$, $\lambda=\tilde{\lambda}$, $\tilde{\lambda}_0=0$, and the
Bautin ideals associated to the holonomies (\ref{bautin4}) are principal.
Consider the circle $S_R$ defined in (\ref{sr}). The Bautin ideal of the
Poincar\'{e} map $P_\lambda$ coincides with the Bautin ideal of the displacement
map $d^1_\lambda-d^2_\lambda$ and
$$
d^1_\lambda(z)-d^2_\lambda(z) =  \varphi(\lambda)( M(z) + R(z,\lambda))
$$
where $\varphi $ is the generator of the Bautin ideal, $R$ is analytic,
$R(z,0)=0$, and $M(z)$ is the Poincar\'{e}-Pontryagin function. As before $M$
satisfies a Fuchs equation with a singularity at $z=0$. We choose $R$ so small,
that the increase of the argument of $M(z)$ along $S_R$, $\arg(z)<\pi$ is close
to the increase of the argument of the leading term of $M$. We note that if
$z^\nu (\log(z))^\mu$ is the leading term of $M$, then
$$
\nu = \nu(P_{\lambda(\varepsilon)})
$$
where $\varepsilon\rightarrow \lambda(\varepsilon)$ is a one-parameter
deformation (a specialization), such that $\varphi(\lambda(.))\neq 0$. We fix
$R$ and choose $\varepsilon_0>0$ so small, that for all $\lambda$, such that
$|\lambda| < \varepsilon_0$, the increase of the argument of the displacement
map $d^1_\lambda(z)-d^2_\lambda(z)$ along $S_R$, $\arg(z)<\pi$, is close to the
increase of the argument of $M(z)$. By making use of Lemma \ref{mainlemma1},
the proof of Theorem \ref{main} is completed as in section~\ref{mainsection}.
Note that the characteristic numbers which appear in the estimate of the
cyclicity may be obtained as characteristic numbers corresponding to a
one-parameter analytic deformation $\mathcal{F}_{\lambda(\varepsilon)}$ of
$\mathcal{F}_0$
$$
\varepsilon\rightarrow \lambda(\varepsilon), \lambda(0)=0, \lim
\tilde{\lambda}(\varepsilon)= \tilde{\lambda}_0
$$
provided that the generators of the Bautin ideals do not vanish identically
along this deformation. Therefore, a multi-parameter version of Theorem
\ref{main} can be formulated as follows
\begin{theorem}
 \label{main+}
 There exists a germ of analytic curve $\varepsilon\rightarrow \lambda(\varepsilon),
 \lambda(0)=0$ in the parameter space, such that
\begin{equation}\label{mainbound+}
Cycl(\Gamma_2,\mathcal{F}_\lambda) \leq 1+ \nu(P_{\lambda(\varepsilon)})+
\max\{\nu(h^{\lambda(\varepsilon)}_{\delta_1}),\nu(h^{\lambda(\varepsilon)}_{\delta_2})\}+
\nu(h^{\lambda(\varepsilon)}_{\delta_1}\circ
h^{\lambda(\varepsilon)}_{\delta_2}) .
\end{equation}
\end{theorem}
\def\cprime{$'$} \def\cprime{$'$} \def\cprime{$'$} \def\cprime{$'$}

\end{document}